\documentclass[11pt]{amsart}
\usepackage{amsmath}
\usepackage{latexsym}
\usepackage{amssymb}
\usepackage{amscd}
\usepackage{amsfonts}
\usepackage[all]{xy}

\DeclareFontShape{OT1}{cmr}{b}{n}{<12> cmr12}{}

\newtheorem{proposition}{Proposition}[section]

\newtheorem{lemma}[proposition]{Lemma}
\newtheorem{definition}[proposition]{Definition}
\newtheorem{lemma-definition}[proposition]{Lemma-Definition}
\newtheorem{theorem}[proposition]{Theorem}

\newtheorem{corollary}[proposition]{Corollary}

\newtheorem{remark}[proposition]{Remark}

\newenvironment{dok}{\par\vspace{-5pt}%
\par\noindent\begingroup%
\leftskip=0em\hspace{0em}{\bf Proof.}}%
{\endgroup\hfill$\Box$}

\newcounter{tmp}

\def\db#1{ \bD^b({#1})}
\def\d#1#2{{\bD}^b_{#1} ({#2})}
\def\perf#1{{\mathfrak P}{\mathfrak e}{\mathfrak r}{\mathfrak f}({#1})}
\def\per#1#2{{\mathfrak P}{\mathfrak e}{\mathfrak r}{\mathfrak f}_{#1}({#2})}
\def\dsing#1{ \bD_{\rm Sg}({#1})}

\def\dhf#1{ {#1}_{hf}}

\def\ove{\overline}

\def\Perf{{\mathfrak P}{\mathfrak e}{\mathfrak r}{\mathfrak f}}

\def\Ho#1,#2,#3,#4{{\operatorname{Hom}}^{#1}_{#2}({#3}\:,\; {#4})}
\def\Ex#1,#2,#3,#4{{\operatorname{Ext}}^{#1}_{#2}({#3}\:,\; {#4})}

\def\lto{\longrightarrow}

\def\A{{\mathcal A}}
\def\B{{\mathcal B}}

\def\D{{\mathcal D}}
\def\F{{\mathcal F}}
\def\E{{\mathcal E}}
\def\G{{\mathcal G}}

\def\K{{\mathcal K}}
\def\N{{\mathcal N}}
\def\O{{\mathcal O}}

\def\J{{\mathcal J}}
\def\T{{\mathcal T}}

\def\P{{\mathcal P}}

\def\P{{\mathcal{P}}}

\def\ZZ{{\mathbb Z}}

\def\bD{{\mathbf D}}
\def\bR{{\mathbf R}}

\def\bL{{\mathbf L}}

\def\ZZ{{\mathbb Z}}
\def\b1{{\mathbf 1}}

\def\bK{{\mathbf H}}
\def\bS{{\mathbf S}}

\def\ZZ{{\mathbb Z}}

\def\PP{{\mathbb P}}

\def\Hom{\operatorname{Hom}}

\def\Ext{\operatorname{Ext}}

\def\Supp{\operatorname{Supp}}

\def\hom{\underline{{\mathcal H}om}}

\def\rI{\operatorname{i}}

\def\Sing{\operatorname{Sing}}

\def\Ker{\operatorname{Ker}\,}

\def\Im{\operatorname{Im}\,}

\def\Spec{\mathbf{Spec}\,}

\def\Inj{\operatorname{Inj}}
\def\inj{\operatorname{inj}}

\def\coh{\operatorname{coh}}
\def\Qcoh{\operatorname{Qcoh}}

\def\wX{\mathfrak{X}}
\def\wJ{\mathfrak{J}}
\def\wF{\mathfrak{F}}
\def\wi{\mathfrak{i}}

\title[]{Formal completions and idempotent completions of triangulated
categories of singularities}

\author[]{Dmitri Orlov}

\address{ Algebra Section, Steklov Mathematical Institute RAN,
Gubkin str. 8,  Moscow 119991, RUSSIA}

\email{orlov@mi.ras.ru}

\thanks{
This work is done under partial financial support by  RFFI grant
08-01-00297, INTAS grant 05-1000008-8118, and  NSh grant
9969.2006.1}

\dedicatory{Dedicated to the blessed memory of my adviser  Vasily Alexeevich
Iskovskikh}

\date{}

\begin{document}

\begin{abstract}
The main goal of this paper is to prove that the idempotent
completions of the triangulated categories of singularities of two
schemes  are equivalent if  the formal completions of these schemes
along singularities are isomorphic. We also discuss Thomason theorem
on dense subcategories and a relation to the negative K-theory.
\end{abstract}
\maketitle
\section{Introduction}

Let $X$ be a noetherian scheme over a field $k.$ Denote by $\db{\coh
X}$
 the bounded derived categories of coherent
 sheaves on $X.$
Since $X$ is noetherian the natural functor from $\db{\coh X}$ to
the unbounded derived category of quasi-coherent sheaves $\bD(\Qcoh
X)$ is fully faithful and realizes  an equivalence of $\db{\coh X}$
with the full subcategory $\bD^{\emptyset,\, b}_{\coh}(\Qcoh X)$
consisting of all cohomologically bounded complexes with coherent cohomologies
(\cite{Il}, Ex.II, 2.2.2).

Denote by $\perf{X}\subseteq \bD(\Qcoh X)$ the full triangulated
subcategory of perfect complexes. Recall that a complex on a scheme
is said to be perfect if it is locally quasi-isomorphic to a bounded
complex of locally free sheaves of finite rank.

The derived category
$\bD(\Qcoh X)$ admits all coproducts and it is well-known that
subcategory of perfect complexes $\perf{X}$ coincides with
subcategory of compact objects in $\bD(\Qcoh X),$ i.e. all objects
$C\in\bD(\Qcoh X)$ for which the functor $\Hom(C,-)$ commutes with
arbitrary coproducts. The category $\perf{X}$ can be considered as a
full triangulated subcategory of $\db{\coh X}.$

\begin{definition}\label{trcsin} We define a triangulated category of singularities of $X,$
denoted by
$\dsing{X},$ as the quotient of the triangulated category $\db{\coh
X}$ by the full triangulated subcategory $\perf{X}.$
\end{definition}

We say that $X$ satisfies a condition {\sf (ELF)} if  $X$ is
 separated  noetherian of finite  Krull dimension and  has enough  locally free sheaves,
i.e.  for any coherent sheaf $\F$ there is an epimorphism
$\E\twoheadrightarrow\F$ with a locally free
 sheaf $\E.$

The last condition implies that any perfect complex
is also globally (not only locally) quasi-isomorphic to a bounded complex of locally free sheaves of finite rank.
For example, any quasi-projective scheme satisfies these conditions.
Note that any closed and any open subscheme of $X$ is also
noetherian, finite dimensional and has enough locally free sheaves.
It is clear for a closed subscheme while for an open subscheme $U$
it follows from the fact that any coherent sheaf on $U$ can be
obtained as the restriction of a coherent sheaf on $X$ (\cite{H}, ex.5.15).

Further in the paper we usually assume that a scheme $X$ satisfies
condition (ELF).

It is  known that if a scheme $X$  is regular  then the
category $\Perf(X)$ coincides with $\db{\coh X}.$ In this case, the triangulated
category of singularities is trivial.

Let $f: X\to Y $ be a morphism of finite Tor-dimension (for example,
a flat morphism or a regular closed embedding). In this case we have an
inverse image functor $\bL f^*:\db{\coh Y}\to \db{\coh X}.$ It is
clear that the functor $\bL f^*$ sends
 perfect complexes on $Y$ to  perfect complexes on $X.$
Therefore, the functor $\bL f^*$ induces  an exact functor $\bL
\bar{f}^* :\dsing{Y} \to \dsing{X}.$

A fundamental  property of  triangulated categories of singularities
is a property of locality in Zarisky topology. It says that for any
open embedding $j: U\hookrightarrow X,$ for which $\Sing(X)\subset
U,$ the functor $\bar{j}^*:\dsing{X}\to \dsing{U}$ is an
equivalence of triangulated categories \cite{Tr}.

On the other hand, two analytically isomorphic singularities can
have non-equivalent triangulated categories of singularities. Even
two different double points given by equations $f=y^2-x^2$ and
$g=y^2-x^2-x^3$ have non-equivalent categories of singularities. The
main reason here is that a triangulated category of singularities
is not necessary idempotent complete. This means
that not for each projector $p: C\to C,\; p^2=p$ there is a decomposition of the form $C=\Ker p\oplus \Im
p.$

For any triangulated category $\T$ we can consider  its so called idempotent completion
 (or Karoubian envelope) $\overline{\T}.$ This
is a category that consists of all kernels of all projectors. It
has a natural structure of a triangulated category and the canonical
functor $\T\to\overline{\T}$ is an exact full embedding \cite{BS}.
Moreover, the category $\overline{\T}$ is idempotent complete now,
i.e. each idempotent $p: C\to C$ in $\overline{\T}$ arises from a
splitting $\Ker p\oplus\Im p.$ We denote by $\overline{\dsing{X}}$
the idempotent completion of the triangulated categories of
singularities.

For any closed subscheme $Z\subset X$ we can define the formal
completion of $X$ along $Z$ as a ringed space $(Z,\,
{\underleftarrow{\lim}} \, \O_{X}/\J^n),$ where $\J$ is
the ideal sheaf corresponding to $Z.$ The formal completion actually
depends only on the closed subset $\Supp Z$ and does not depend on a
scheme structure on $Z.$ We denote by $\wX$ the formal completion of
$X$ along its singularities $\Sing(X).$

The main goal of this paper is to prove that for any two schemes $X$ and $X'$ satisfying (ELF), if the formal
completions $\wX$ and $\wX'$ along singularities are isomorphic,
then the idempotent completions of the triangulated categories of
singularities $\overline{\dsing{X}}$ and $\overline{\dsing{X'}}$ are
equivalent (Theorem \ref{main}). Actually, we show a little bit more. We prove that
any object of $\dsing{X}$ is a direct summand of an object from its
full subcategory $\d{\Sing(X)}{\coh X}/\per{\Sing(X)}{X},$ where $\d{\Sing(X)}{\coh X}$ and
$\per{\Sing(X)}{X}$ are subcategories of $\db{\coh X}$ and $\Perf(X)$ respectively, consisting of complexes with cohomologies supported on
$\Sing X$ (Proposition \ref{idemp}).

Thus, to any scheme $X$ we can attach the category
$\overline{\dsing{X}}$ and two subgroups in $K_0 (\overline{\dsing{X}})$ that
by Thomason theorem \cite{T} (see Theorem \ref{Thomason}) one-to-one corresponds to
the dense subcategories $\dsing{X}$ and $\d{\Sing(X)}{\coh X}/\per{\Sing(X)}{X}$ respectively.
We discuss this correspondence and a relation to the negative K-theory of the category of perfect complexes in  the last section.

\section{Completions}
Let $X$ be a noetherian scheme and let $i: Z {\hookrightarrow} X$ be
a closed subspace. Let $\coh_Z X\subset \coh X$ be the abelian
subcategory of coherent sheaves on $X$ with support on $Z.$

Consider the natural functor from $\db{\coh_Z X}$ to $\db{\coh X}.$
It can be easily shown that this functor is fully faithful and gives
an equivalence with the full subcategory $\d{Z}{\coh X}\subset
\db{\coh X}$ consisting of all complexes cohomologies of which are
supported on $Z.$ (In other words, the subcategory $\d{Z}{\coh X}$
consists of all complexes restriction of which on the open subset
$U=X\backslash Z$ is acyclic.)

At first, let us consider abelian category of quasi-coherent sheaves
$\Qcoh X$ and its abelian subcategory $\Qcoh_Z X$ of quasi-coherent
sheaves with support on $Z$ (or $Z$-torsion sheaves), i.e. all
quasi-coherent sheaves $\F$ such that $j^*\F=0,$ where $j: U\to X$
is the open embedding of the complement $U=X\backslash Z.$ The
inclusion functor $\rI: \Qcoh_Z X\to \Qcoh X$ has a right adjoint
$\varGamma_Z$ which associates to each quasi-coherent sheaf $\F$ its
subsheaf of sections with support in $Z.$ It can be shown that for a
quasi-coherent sheaf $\F$ we have an isomorphism
$$
\varGamma_Z
(\F)\cong\underset{n}{\underrightarrow{\lim}}\;\hom_{\O_X}(\O_X/\J^n,\,
\F),
$$
where $\J$ is a some ideal sheaf such that $Z=\Supp(\O_X/\J).$ The
functor $\varGamma_Z$ has a right-derived functor $\bR\varGamma_Z :
\bD(\Qcoh X)\to \bD(\Qcoh_Z X)$  via h-injective resolutions
\cite{Sp}.

It is known that the canonical functor $\rI: \bD(\Qcoh_Z
X)\to\bD(\Qcoh X)$ is fully faithful and realizes equivalences of
$\bD(\Qcoh_Z X)$ with the full subcategory $\bD_{Z}(\Qcoh X)$
consisting of all complexes cohomologies of which are supported on
$Z.$ It is proved for noetherian schemes for example in \cite{AJL1}
(Prop. 5.2.1 and Prop. 5.3.1). (It is also true for quasi-compact
and separated X and proregular embedded $Z\subset X$ as shown in
\cite{AJL2}.) To prove this fact it is sufficient to show that for
any $C^{\cdot}\in\bD_{Z}(\Qcoh X)$ the natural map
$\rI\bR\varGamma_Z (C^{\cdot}) \to C^{\cdot}$ is an isomorphism
(see, for example, \cite{AJL1} Lemma 5.2.2). Since the functor
$\bR\varGamma_Z$ is bounded for noetherian schemes by usual "way
out" argument (\cite{Ha}, \S 7)  it is sufficient to check the
isomorphism $\rI\bR\varGamma_Z (\F) \to \F$ only for sheaves $\F\in
\Qcoh_Z X.$ That is evident, because $\bR\varGamma_Z \F\cong
\varGamma_Z \F,$ when $\F$ is $Z$-torsion.

Thus, for any object $C^{\cdot}\in \bD(\Qcoh X)$ there is a
distinguished triangle of the form
$$
\rI\bR\varGamma C^{\cdot}\to C^{\cdot}\to \bR j_*j^* C^{\cdot},
$$
which shows that the categories $\bD(\Qcoh U)$ and $\bD(\Qcoh_Z X)$
are equivalent to the quotient categories $\bD(\Qcoh X)/\bD_Z(\Qcoh
X)$ and $\bD(\Qcoh X)/\bD(\Qcoh U)$ respectively.

\begin{lemma}\label{support} Let $X$ be a noetherian scheme and let $Z$ be a closed
subspace. Then the natural functor $\db{\coh_Z X}\to\db{\coh X}$ is
fully faithful and gives an equivalence with the full
subcategory $\d{Z}{\coh X}\subset \db{\coh X}$ consisting of all
complexes cohomologies of which are supported on $Z.$
\end{lemma}
\begin{dok}
We know that the natural functors $\db{\coh X}\hookrightarrow
\bD(\Qcoh X)$ and $\bD(\Qcoh_Z X)\hookrightarrow \bD(\Qcoh X)$ are
fully faithful. This implies that the functor $\db{\coh_Z
X}\to\db{\coh X}$ is fully faithful iff the functor $\db{\coh_Z
X}\to \bD(\Qcoh_Z X)$ is fully faithful. Denote by $\phi$ the
natural embedding of $\coh_Z X$ to $\Qcoh_Z X.$ Since coherent
sheaves generate the category $\db{\coh_Z X}$ it is enough to show
that for any two coherent sheaves $\F, \G\in \coh_Z X$ the natural
maps $\Ext^n (\F, \G)\to\Ext^n (\phi(\F), \phi(\G))$ are
isomorphisms. We know that it is evidently true for $n=0.$ Now to
apply induction, it is sufficient to check that for any $e\in \Ext^n
(\phi(\F), \phi(\G)), \; n\ge 1,$ there is an epimorphism $\F'\to
\F$ which erases $e$ (\cite{Il}, Ex. II, Lemma 2.1.3). Any such
element $e$ can be represented by an exact sequence in $\Qcoh_Z X$
$$
0\to\G\to \E_{n-1}\to\cdots\to\E_0\to\F\to 0,
$$
where $\E_i$ are quasi-coherent sheaves with support on $Z.$ The
epimorphism $\E_0\to \F$ erases $e.$ Since any quasi-coherent sheaf
on a noetherian scheme is a direct colimit of its coherent
subsheaves there is a coherent subsheaf $\F'\subset\E_0$ that covers
$\F.$ As a subsheaf of $\E_0$ it is erase $e$ and   belongs to
$\coh_Z X.$ Thus, the natural functor $\db{\coh_Z X}\to\db{\coh X}$
is fully faithful and gives an equivalence with the full subcategory
$\d{Z}{\coh X}\subset \db{\coh X}$ consisting of all complexes
cohomologies of which are supported on $Z.$
\end{dok}

The restriction functor $j^*$ sends coherent sheaves to coherent and
we get a functor from the quotient category $\bD^b(\coh
X)/\bD^b_Z(\coh X)$ to the derived category $\bD^b(\coh U).$ This functor
establishes an equivalence between these categories.

\begin{lemma}
 Let $X$ be a noetherian scheme. Then the natural functor
$
\bD^b(\coh X)/\bD^b_Z(\coh X)\to\bD^b(\coh U)
$
is an equivalence.
\end{lemma}
\begin{dok} This fact is known and we omit a proof. There are few different ways to get it.

First, we know this fact for quasi-coherent sheaves and by Lemma \ref{ffemb} it is enough to show that
any map from a bounded complex of coherent sheaves to an object of $\bD_Z(\Qcoh X)$ admits a factorization
through an object from $\bD^b_Z(\coh X).$

Second, since $\coh U$ is the quotient of the abelian category $\coh X$ by the Serre subcategory $\coh_Z X$ it is possible to show
that the functor $\bD^b(\coh X)/\bD^b_Z(\coh X)\to\bD^b(\coh U)$ is surjective on objects and on morphisms.
This implies an equivalence as well.
\end{dok}

\begin{remark}{\rm
Note that the second way allows us to prove a more general result.
For any Serre subcategory $\B$ of an abelian category $\A$ the functor $F: \db{\A}/\bD^b_{\B}(\A)\to
\db{\A/\B}$ is an  equivalence of triangulated categories, where $\bD^b_{\B}(\A)$  is the
full subcategory in $\db{\A}$ consisting of all complexes with cohomologies
in ${\B}$ (it is known but unpublished \cite{BO}).}
\end{remark}

Denote by $\per{Z}{X}$  the intersection $\perf{X}\cap\d{Z}{\coh
X}.$

\begin{lemma}\label{perfc} Let $X$ satisfies (ELF). Then an object $A\in\d{Z}{\coh X}$
belongs to $\per{Z}{X}$  iff  for any object
$B\in\d{Z}{\coh X}$ all $\Hom(A, B[i])$ are trivial except for
finite number of $i\in\ZZ.$
\end{lemma}
\begin{dok}
Denote by $\dhf{\D}\subset \d{Z}{\coh X}$ the full subcategory
consisting of all objects $A$ such that for any object
$B\in\d{Z}{\coh X}$ the spaces $\Hom(A, B[i])$ are trivial except
for finite number of $i\in\ZZ.$ If an object $A\in\per{Z}{X}$  then
it is quasi-isomorphic to a bounded complex of vector bundles. Since
the cohomologies of any coherent sheaf is bounded by the Krull
dimension of the scheme we have that for any vector bundle $\P$ and
any coherent sheaf $\F$ there is an equality $\Ext^i(\P, \F)=0$ when
$i$ is greater than Krull dimension of $X.$ Therefore, $A$ belongs
to the subcategory $\dhf{\D}.$

Suppose now that $A\in \dhf{\D}.$ The object $A$ is a bounded
complex of coherent sheaves. Let us take locally free  bounded above
resolution $P^{\cdot}\stackrel{\sim}{\to} A$ and consider a good
truncation $\tau^{\ge -k}P^{\cdot}$ for sufficient large $k\gg 0$
which is clearly isomorphic to $A$ in $\D.$

Since $A\in \dhf{\D},$ for any closed point $t: x\hookrightarrow X$
the groups $\Hom(A, t_*\O_x[i])$ are zero for $|i|\gg 0.$ This means
that for sufficiently large $k\gg 0$ the truncation $\tau^{\ge
-k}P^{\cdot}$ is a complex of locally free sheaves at the point $x,$
and, hence, in some neighborhood of $x.$ The scheme $X$ is
quasi-compact. This implies that there is a common sufficiently
large $k$ such that the truncation $\tau^{\ge -k}P^{\cdot}$ is a
complex of locally free sheaves everywhere on $X,$ i.e. $A$ is
perfect.
\end{dok}

 The natural embedding $\d{Z}{\coh X}\hookrightarrow \db{\coh X}$
induces a functor between quotient categories
$$
\d{Z}{\coh X}/\per{Z}{X}\lto\dsing{X}:=\db{\coh X}/\perf{X}.
$$
It can be proved that this functor between quotient categories is
fully faithful too.
To prove it we need the
following well-known lemma.

\begin{lemma}{\rm (\cite{Ve, KS})}
\label{ffemb} Let $\D$ be a triangulated category and $\D', \N$ be
full triangulated subcategories. Let $\N'=\D'\cap \N.$ Assume that
any morphism $N\to X'$ (resp. any morphism $X'\to N$) with $N\in \N$
and $X'\in D'$ admits a factorization $N\to N'\to X'$ (resp. $X'\to
N'\to N$) with $N'\in \N'.$ Then the natural functor
$
\D'/\N'\lto \D/\N
$
is fully faithful.
\end{lemma}

\begin{lemma}\label{fulf}The functor $\d{Z}{\coh X}/\per{Z}{X}\to\dsing{X}$ is  fully faithful.
\end{lemma}
\begin{dok}
By Lemma \ref{ffemb} we should show that any morphism $P^{\cdot}\to
C^{\cdot},$ where $P^{\cdot}$ is a perfect complex and $C^{\cdot}$
is an object of $\d{Z}{\coh X},$ can be factorized through an object
${P^{\cdot}}'\in\per{Z}{X}.$ Since there is an equivalence
$\db{\coh_Z X}\cong\d{Z}{\coh X}$ we can assume that the object
$C^{\cdot}$ is a bounded complex of coherent sheaves with support on
$Z.$ This implies that there is a subscheme structure $i_S:
S\hookrightarrow X$ with support on $Z$ such that $C^{\cdot}\cong
\bR i_{S*} C^{\cdot'}.$ Consider a vector bundle $\E$ on $X$ that
covers the ideal sheaf $\J_S.$ It exists by (ELF) condition. Denote
by $K^{\cdot}$ the Koszul complex
$0\to\det(\E)\to\cdots\to\E\to\O_X\to 0.$ We have canonical maps
$P^{\cdot}\to P^{\cdot}\otimes \K^{\cdot}\to P^{\cdot}\otimes O_S.$
Now any map from $P^{\cdot}$ to $C^{\cdot}\cong \bR i_{S*}
C^{\cdot'}$ is factorized through $P^{\cdot}\otimes O_S$ and, hence,
through $P^{\cdot}\otimes \K^{\cdot}.$ But the object
$P^{\cdot}\otimes \K^{\cdot}$ is perfect as a tensor product of two
perfect complexes and has cohomologies with supports on $Z.$ Thus,
it belongs to $\per{Z}{X}.$
\end{dok}

Now consider the case when $Z$ is exactly the subset of
singularities of $X,$ i.e. $Z=\Sing(X).$

\begin{proposition}\label{idemp}
Any object of $\dsing{X}$ is a direct summand of an object from its
full subcategory $\d{\Sing(X)}{\coh X}/\per{\Sing(X)}{X}.$ In
particular, idempotent completions of these categories are
equivalent.
\end{proposition}
\begin{dok}
Since any object of the category $\dsing{X}$ is represented by a
coherent sheaf up to shift (\cite{Tr}, Lemma 1.11) it is sufficient
to consider the coherent sheaf $\F.$ Let us take the locally free
resolution $P^{\cdot}\to \F$ and consider the brutal truncation
$\sigma^{\ge-n}P^{\cdot}$ for sufficiently large $n>\dim X.$ Denote
by $\G$ the (-n)-th cohomology of $\sigma^{\ge-n}P^{\cdot}.$ Let $\alpha:\F\to\G[n+1]$
be the corresponding map in $\db{\coh X}.$ Its image
in the category $\dsing{X}$ is an isomorphism. On the other hand,
consider the functor $j^*: \db{\coh X}\to \db{\coh U},$ where
$U=X\backslash \Sing(X).$ Since $U$ is smooth and $n>\dim U$ the
image $j^{*}(\alpha)$ is zero. But the category $\db{\coh U}$ is the
quotient of the category of $\db{\coh X}$ by the subcategory
$\d{\Sing(X)}{\coh X}.$ Hence the morphism $\alpha$ is factorized
through an object $A$ of $\d{\Sing(X)}{\coh X}.$ Therefore, in the
quotient category $\dsing{X}$ the object $\F$ is a direct summand of
the image of the object $A$ in $\dsing{X}.$
\end{dok}

For any scheme $X$ we denote by $\wX_Z$ the formal completion of $X$
along a closed subspace $Z$ and denote by  $\kappa: \wX_Z\to X$ the
canonical morphism. Let $\J$ be an ideal sheaf such that
$\Supp(\O_X/\J)=Z$ and let $\wJ$ be a corresponding ideal of
definition of the formal noetherian scheme $\wX_Z.$ We set
$$
\varGamma_{\wX}
(\wF):=\underset{n}{\underrightarrow{\lim}}\;\hom_{\O_{\wX_Z}}(\O_{\wX_Z}/\wJ^n,\,
\wF),
$$
for any quasi-coherent sheaf $\wF$ on $\wX_Z.$ This functor depends
only on $\O_{\wX_Z}$ and does not depend on the ideal $\wJ.$ We say
that $\wF\in \Qcoh \wX_Z$ is a torsion sheaf if
$\varGamma_{\wX}(\wF)=\wF.$ We denote by $\coh_t \wX_Z$ (resp.
$\Qcoh_t \wX_Z$) the full subcategory of $\Qcoh \wX_Z$  whose
objects are the (quasi)-coherent torsion sheaves. It is easy to see
that under the inverse image functor $\kappa^*$ a $Z$-torsion
(quasi)-coherent sheaf on $X$ goes to a torsion (quasi)-coherent
sheaf on $\wX_Z.$ Indeed, applying $\underrightarrow{\lim}$ to the
isomorphisms
$$
\kappa^* \hom_X(\O_X/\J^n, \F)\stackrel{\sim}{\lto}
\hom_{\wX_Z}(\O_{\wX}/{\wJ}^n, \kappa^*\F)
$$
we get a natural isomorphism
$\kappa^*\varGamma_Z\cong\varGamma_{\wX}\kappa^*.$ Hence, we obtain
that $\kappa^*(\Qcoh_Z X)\subset\Qcoh_t \wX_Z$ and $\kappa^*(\coh_Z
X)\subset\coh_t \wX_Z.$

Now, if $\F$ is a $Z$-torsion coherent sheaf on $X$ then there is an
integer $n$ such that $\F$ comes from $Z_n=\Spec \O_X/\J^n=(\wX_Z,
\O_{\wX_Z}/\wJ^n)$ under the closed inclusion $i_n: Z_n\to X,$ i.e.
$\F=i_{n*}\F'$ for some coherent sheaf $\F'\in\coh Z_n.$ Consider
the cartesian diagram
$$
\begin{CD}
Z_n@>\mathfrak{i}_n>>& \wX_Z\\
@|&@VV{\kappa}V\\
Z_n@>i_n>>& X
\end{CD}
$$

We have a sequence of isomorphisms
$
\kappa_*\kappa^*\F\cong\kappa_*\kappa^*i_{n*}\F'\cong
\kappa_*\wi_{n*}\F'\cong i_{n*}\F'\cong \F.
$
If now $\wF$ is a
torsion sheaf on $\wX_Z$ then again there is an integer $n$ such
that $\wF\cong \wi_{n*}\F'$ and
$\kappa^*\kappa_*\wF\cong\kappa^*\kappa_*\wi_{n*}\F'\cong
\kappa^*i_{n*}\F'\cong \wi_{n*}\F'\cong \wF.$
Thus, we obtain that
the functors $\kappa^*$ and $\kappa_*$ induce inverse equivalences
between the abelian categories $\coh_Z X$ and $\coh_t \wX_Z.$ It is
also can be shown that the functor $\kappa_*$ sends quasi-coherent
torsion sheaf to quasi-coherent $Z$-torsion sheaves, because
$\kappa_*$ commutes with colimits (see \cite{AJL1}, Prop. 5.1.1,
5.1.2). Thus, we get the following proposition.

\begin{proposition}{\rm (\cite{AJL1})}\label{ajl}
Let $X$ be a noetherian scheme and $\kappa: \wX_Z\to X$ be a formal
completion of $X$ along a closed subspace $Z.$ Then the functors
$\kappa^*$ and $\kappa_*$ induce  inverse equivalences between the
categories $\coh_Z X$ and $\coh_t \wX_Z,$ and between the categories
$\Qcoh_Z X$ and $\Qcoh_t \wX_Z.$
\end{proposition}

\begin{corollary}\label{eqsup} Let $X$ and $X'$ be two schemes satisfying (ELF). Assume that the formal schemes
$\wX_{Z}$ and $\wX'_{Z'}$ are isomorphic. Then the derived
categories $\d{Z}{\coh X}$ and $\d{Z'}{\coh X'}$ (resp. $\d{Z}{\Qcoh
X}$ and $\d{Z'}{\Qcoh X'}$) are equivalent.
\end{corollary}
\begin{dok} By Lemma \ref{support} there is an equivalence $\d{Z}{\coh X}\cong\db{\coh_Z X}$ (resp. $\d{Z}{\Qcoh X}\cong\db{\Qcoh_Z X}$
by  and by Proposition \ref{ajl} we have $\coh_Z X\cong\coh_t \wX_Z$
(resp. $\Qcoh_Z X \cong \Qcoh_t \wX_Z$). Since $\wX_Z\cong\wX'_{Z'}$
we obtain that $\coh_Z X\cong\coh_{Z'} X'$ (resp. $\Qcoh_Z
X\cong\Qcoh_{Z'} X'$). Therefore, the derived categories are
equivalent as well.
\end{dok}

\begin{theorem}\label{main} Let $X$ and $X'$ be two schemes satisfying (ELF). Assume that the formal completions
$\wX$ and $\wX'$ along singularities are isomorphic. Then the
idempotent completions of the triangulated categories of
singularities $\overline{\dsing{X}}$ and $\overline{\dsing{X'}}$ are
equivalent.
\end{theorem}
\begin{dok}
By Corollary \ref{eqsup} there is an equivalence between
$\d{\Sing(X)}{\coh X}$ and $\d{\Sing(X')}{\coh X'}.$ By Lemma
\ref{perfc} the subcategories $\Perf(X)$ and $\Perf(X')$ are also
equivalent, because they can be defined in the internal terms of the
bounded derived categories of coherent sheaves with support on $Z$ and $Z'.$ Hence,
there is an equivalence between quotient categories
$$
\d{\Sing(X)}{\coh
X}/\per{\Sing(X)}{X}\stackrel{\sim}{\lto}\d{\Sing(X')}{\coh
X'}/\per{\Sing(X')}{X'},
$$
It induces an equivalence between their idempotent completions, which
by Proposition \ref{idemp}
coincide with the idempotent completions of the triangulated categories of
singularities.
\end{dok}

\section{Localization in Nisnevich topology and isomorphisms infinitely  near singularities}

Let $X$ be a noetherin scheme and $i: Z {\hookrightarrow} X$ be a
closed subscheme. Consider the pair $(Z, X).$ Let $f: X'\to X$ be a
map of schemes.

\begin{definition} We say that $f$ is an isomorphism
infinitely near $Z$ if  it is flat over $Z$ and the fiber  product
$Z'=Z\times_X X'$ is isomorphic to $Z.$
\end{definition}

It can be proved that this condition does not depend on a choice of
a closed subscheme with the underline subspace $\Supp Z$
(\cite{TT}, Lemma 2.6.2.2). In particular, if it holds for $Z=\Spec
\O_X/\J$ it also holds for infinitesimal thickenings $Z_n:=\Spec
\O_X/\J^n.$ Thus, we may say that $f$ is an isomorphism infinitely
near the closed subspace $\Supp Z.$

This implies that any morphism $f: X'\to X,$ which is an isomorphism
infinitely near $Z,$ induces an isomorphism between the formal
completions $\mathfrak{f}:
\mathfrak{X}'_{Z'}\stackrel{\sim}{\to}\mathfrak{X}_Z.$ Hence, by
Corollary \ref{eqsup} we obtain that the derived categories
$\d{Z}{\coh X}$ and $\d{Z'}{\coh X'}$ are equivalent for any
morphism $f: X'\to X$ that is an isomorphism  infinitely near $Z$
(see \cite{TT}, Th. 2.6.3).

Important examples of such morphisms are Nisnevich neighborhoods of
$Z$ in $X.$
\begin{definition} An $X$\!-scheme $\pi: Y\to X$ is called a Nisnevich neighborhood
of $Z$ in $X$ if the morphism $\pi$ is etale and the fiber product
$Z\times_X Y$ is isomorphic to $Z.$
\end{definition}

\begin{proposition}\label{qcoheq}
Let a scheme $X$ satisfy (ELF) and let $Z\subset X$ be a closed
subscheme. Then for any morphism $f: X'\to X,$ that is an
isomorphism infinitely near of $Z,$ the functors $f^*: \bD_Z(\Qcoh
X)\to\bD_Z(\Qcoh X')$ and $f^*: \bD^b_Z(\coh X)\to\bD^b_Z(\coh X')$
are equivalences.
\end{proposition}
\begin{dok}
The morphism $f: X'\to X$ induces an isomorphism between the formal
completions $\mathfrak{f}:
\mathfrak{X}'_{Z'}\stackrel{\sim}{\to}\mathfrak{X}_Z$
(\cite{TT}, Lemma 2.6.2.2). Hence, by Corollary \ref{eqsup} we obtain
that the derived categories $\d{Z}{\coh X}$ and $\d{Z}{\coh X'}$
(resp. $\bD_Z(\Qcoh X)$ and $\bD_Z(\Qcoh X')$) are equivalent.
\end{dok}
\begin{proposition}
Let schemes $X$ and $X'$ satisfy (ELF) and let $f: X'\to X$ be a
morphism that is an isomorphism  infinitely near $Z.$ Suppose the
complements $X\backslash Z$ and $X'\backslash Z$ are smooth. Then
the functor $\bar{f}^*: \dsing{X}\to\dsing{X'}$ is fully faithful
and, moreover, any object $B\in \dsing{X'}$ is a direct summand of a
some object of the form $\bar{f}^* A.$
\end{proposition}
\begin{dok}
By assumption $\Sing(X)\cong\Sing(X')\subseteq Z,$ hence, $f$ is an
isomorphism infinitely near $\Sing(X).$ By Proposition \ref{qcoheq}
and Lemma \ref{perfc} we obtain an equivalence
$$
\d{\Sing(X)}{\coh
X}/\per{\Sing(X)}{X}\stackrel{\sim}{\lto}\d{\Sing(X')}{\coh
X'}/\per{\Sing(X')}{X'}.
$$
By Proposition \ref{idemp} their idempotent completions coincide
with idempotent completions of triangulated categories of
singularities. Hence, the natural functor $\bar{f}^*$ is fully
faithful and any object $\dsing{X'}$ is a direct summand of an
object from $\dsing{X}.$
\end{dok}
\begin{corollary}
Let a scheme $X$ satisfy (ELF) and the complement $X\backslash Z$ is
smooth. Then for any Nisnevich neighborhood $\pi: Y\to X$  of $Z$
the functor $\bar{\pi}^*: \dsing{X}\to\dsing{Y}$  is fully faithful
and, moreover, any object $B\in \dsing{Y}$ is a direct summand of a
some object of the form $\bar{\pi}^*A.$
\end{corollary}

\begin{remark}{\rm
Let $(A, p)$ be a pair consisting of a commutative $k$\!-algebra of finite
type $A$ and a prime ideal $p.$
Consider the henselization $(A_h, p_h)$ of this pair. By definition,
$A_h=\underrightarrow{\lim} B,\; p_h=p A_h,$
where the limit is taking by the category of all Nisnevich neighborhoods
$\Spec B\to \Spec A$ of $\Spec A/p$ in $\Spec A.$
In particular, we have $A/p\stackrel{\sim}{\to} A_h/p_h.$
Let $\hat{A}$ be the $p$\!-adic completion of $A.$ By one of application
of Artin Approximation (Theorem 3.10, \cite{Ar})
for any finitely generated $\hat{A}$\!-module $\bar{M},$
which is locally free on $\Spec \hat{A}$ outside $V(\hat{p}),$
there is an $A_h$\!-module
$M$ such that $\hat{M}\cong \bar{M}.$ Assume now that $\Spec A$
is regular outside $V(p).$ This implies that
the natural functor from  $\dsing{\Spec A_h}$ to $\dsing{\Spec \hat{A}}$
is an equivalence, because any object
of a triangulated category of singularities can be represented by a coherent sheaf
which is locally free on the complement to the singularities.
If moreover, $K_{-1}(A_h)=0$ then the triangulated category $\dsing{\Spec A_h}$
is idempotent complete (see next section), i.e. it coincides with
$\overline{\dsing{\Spec A}}.$ For example, it is true
when $A$ is a normal local ring of dimension two of essentially finite type \cite{W}. }
\end{remark}

\section{Thomason theorem and groups $K_{-1}$}
A full triangulated subcategory $\N$ of a triangulated category $\T$
is called {\sf dense} in $\T$ if each object of $\T$ is a direct
summand of an object isomorphic to an object in $\N.$ There is a not
so well-known but amazing theorem of R.~Thomason which allows us to
describe all strictly full
dense subcategories in a triangulated category.

\begin{theorem}\label{Thomason}{\rm (R.~Thomason, \cite{T})}
Let $\T$ be an essentially small triangulated category. Then there
is a one-to-one correspondence between the strictly full dense
triangulated subcategories $\N$ in $\T$ and the subgroups $H$ of the
Grothendieck group $K_0(\T).$

To $\N$ corresponds the subgroup which is the image of $K_0(\N)$ in
$K_0(\T).$ To $H$ corresponds the full subcategory $\N_H$ whose
objects are those $N$ in $\T$ such that $[N]\in H \subset K_0(\T).$
\end{theorem}
\begin{remark}{\rm
Recall that a full triangulated subcategory $\N$ of $\T$ is called
{\sf strictly} full if it contains every object of $\T$ that is isomorphic
to an object of $\N.$}
\end{remark}

Thus, to any scheme $X$ we can attach  the triangulated category $\overline{\dsing{X}}$
and two subgroups in the Grothendieck group $K_0(\overline{\dsing{X}})$ which are related to
the natural dense subcategories $\dsing{X}$ and $\d{\Sing(X)}{\coh
X}/\per{\Sing(X)}{X}$ and which by Thomason's theorem uniquely determine them.

The sequence of triangulated categories
\begin{equation}\label{exseq}
\Perf(X)\lto\db{\coh X}\lto \overline{\dsing{X}}
\end{equation}
is exact in Definition 1.1 of \cite{Schl}, i.e. the first functor is a full embedding and
the quotient of this map is dense subcategory in the third category.

Following Amnon Neeman \cite{Ne} this exact sequence can be
considered as an exact sequence of triangulated categories of
compact objects coming from a localizing sequence of compactly
generated triangulated categories. As we know the category of
perfect complexes $\Perf(X)$ is the category of compact objects in
$\bD(\Qcoh X).$ It is proved by H.~Krause \cite{Kr} that the category
$\db{\coh X}$ can be considered as the category of compact object in
the homotopy category of injective quasi-coherent sheaves
$\bK(\Inj X)$ for a noetherian scheme X. On the other hand, the derived category $\bD(\Qcoh X)$
is equivalent to the full subcategory $\bK_{\inj}(\Qcoh X)\subset \bK(\Inj X)$
of h-injective complexes, i.e. such complexes $I$ that $\Hom_{\bK(\Qcoh X)}(A, I)=0$ for all acyclic complexes
$A$ from homotopy category $\bK(\Qcoh X)$ (see \cite{Sp}).
Since $X$ is noetherian the category $\bK_{\inj}(\Qcoh X)$ closed with respect to formation of coproducts
 and, furthermore,
the inclusion functor $\bK_{\inj}(\Qcoh X)\hookrightarrow \bK(\Inj X)$ respects coproducts.
This means that $\bK_{\inj}(\Qcoh X)$ is localizing subcategory of $\bK(\Inj X)$ and we have a localizing sequence
$$
\bK_{\inj}(\Qcoh X)\lto \bK(\Inj X)\stackrel{Q}{\lto} \bK(\Inj X)/\bK_{\inj}(\Qcoh X).
$$
Moreover, the quotient functor $Q$ has a right adjoint, which is called Bousfield localizing functor.
It identifies the quotient category
$\bK(\Inj X)/\bK_{\inj}(\Qcoh X)$ with the triangulated category of all acyclic complexes of injective objects $\Inj X\subset \Qcoh X.$
The latter category is called stable derived category and will be denoted by $\bS(\Qcoh X).$
By Theorem 2.1 of \cite{Ne} the idempotent completion $\overline{\dsing{X}}$ is equivalent
to the category of all compact objects in the stable derived category $\bS(\Qcoh X)$ (for more details see \cite{Kr}).

By Theorem 11.10 of
\cite{Schl} the sequence (\ref{exseq}) induces a long exact sequence for K-groups
$$
K_0(\Perf(X))\lto K_0(\db{\coh X})\lto K_0(\overline{\dsing{X}})\lto
K_{-1}(\Perf(X))\lto 0.
$$
Here we used  Theorem 9.1 from \cite{Schl} asserting that $K_{-1}$
for a small abelian category is trivial. Therefore, we obtain a
short exact sequence
$$
0\lto K_0(\dsing{X})\lto K_0(\overline{\dsing{X}})\lto
K_{-1}(\Perf(X))\lto 0,
$$
which shows that $K_{-1}(\Perf(X))$ is a measure of the difference
between $\dsing{X}$ and its idempotent
completion $\overline{\dsing{X}}.$ By the same reason, we have another short exact sequence
$$
0\lto K_0(\d{\Sing(X)}{\coh X}/\per{\Sing(X)}{X})\lto
K_0(\overline{\dsing{X}})\lto K_{-1}(\Perf_{\Sing(X)}(X))\lto 0.
$$
Now a long exact sequence for $U=X\backslash \Sing(X)$
$$
K_0(\Perf(X))\to K_0(\Perf(U))\to K_{-1}(\Perf_{\Sing(X)}(X))\to
K_{-1}(\Perf(X))\to 0,
$$
which follows from the Thomason's Localization Theorem 7.4
\cite{TT}, shows a difference between $K_{-1}(\Perf_{\Sing(X)}(X))$
and $K_{-1}(\Perf(X)).$

The negative K-groups (which is due to Bass) are defined from the
following exact sequences
\begin{multline*} 0\to K_i(\Perf(X))\to K_i(\Perf(X[t]))\oplus
K_i(\Perf(X[t^{-1}]))\to K_i(\Perf(X[t,t^{-1}]))\to\\
\to K_{i-1}(\Perf(X))\to 0.
\end{multline*}

In particular, the group  $K_{-1}(\Perf(X))$ is isomorphic to the
cokernel of the canonical map $K_0(\Perf(X[t]))\oplus
K_0(\Perf(X[t^{-1}]))\to K_0(\Perf(X[t,t^{-1}])).$

By Theorem \ref{main} we know that for any two schemes $X$ and $X',$
the formal completions of which along singularities are isomorphic, we have
$$
\overline{\dsing{X}}\cong \overline{\dsing{X'}}\quad \text{and}
\quad \d{\Sing(X)}{\coh X}/\per{\Sing(X)}{X}\cong \d{\Sing(X')}{\coh
X'}/\per{\Sing(X')}{X'}.
$$
On the other hand, in this case the triangulated categories of
singularities $\dsing{X}$ and $\dsing{X'}$ are not necessary
equivalent as we know.

There is also another type of relations between schemes which give equivalences
for triangulated categories of singularities but under which the quotient  categories
$\d{\Sing(X)}{\coh X}/\per{\Sing(X)}{X}$ are not necessary equivalent. It is described in \cite{Or}.

Let $S$ be a  noetherian regular scheme. Let $\E$ be a vector bundle
on $S$ of rank $r$ and let $s \in H^0(S, \E)$ be a  section. Denote
by $X\subset S$ the zero subscheme of $s.$ Assume that the section
$s$ is regular, i.e. the codimension of the subscheme $X$ in $S$ coincides
with the rank $r.$

Consider the  projective bundles $S'=\PP(\E^{\vee})$
and $T=\PP(\E^{\vee}|_X),$ where $\E^{\vee}$ is the dual
bundle.

The section $s$ induces a section $s'\in H^0( S', \O_{\E}(1))$ of
the Grothendieck line bundle $\O_{\E}(1)$ on $S'.$ Denote by $Y$ the
divisor on $S'$ defined by the section $s'.$ The natural closed
embedding of $T$ into $S'$ goes through $Y.$
All
schemes defined above can be included in the following commutative
diagram.

$$
\xymatrix{ T
\ar[d]_{p} \ar@{->}[r]^{i}
& Y
 \ar[rd]^{\pi}\ar@{->}[r]^{u} & S'
 \ar[d]^{q}
\\
X \ar@{->}[rr]^{j} &  & S}
$$

Consider the composition functor $\bR i_* p^*: \db{\coh X}\to
\db{\coh Y}$ and denote it by $\Phi_{T}.$
\begin{theorem}{\rm (\cite{Or})}\label{main1}
Let schemes $X, Y,$ and $T$ be as above. Then the functor
$$\Phi_{T}:\db{\coh(X)}\lto \db{\coh(Y)}$$ defined by the formula $
\Phi_{T}(\cdot)=\bR i_* p^*(\cdot) $ induces a functor
$$
\ove{\Phi}_{T}: \dsing{X}\lto \dsing{Y},
$$
which is an equivalence of triangulated categories.
\end{theorem}
The functor $\Phi_{T}=\bR i_{*} p^*$ has a right adjoint functor
which we denote by $\Phi_{T *}.$ It can be represented as a
composition $\bR p_{*} i^{\flat},$ where $i^{\flat}$ is right
adjoint to $\bR i_{*}.$ Functor $i^{\flat}$ has the form $\bL
i^*(\cdot )\otimes \omega_{T/Y}[-r+1],$ where $\omega_{T/Y}\cong
\Lambda^{r-1}\N_{T/Y}$ is the relative dualizing sheaf.

It is easy to see that all singularities of $Y$ are concentrated  over the singularities of $X,$ hence the functor $\Phi_{T
*}=\bR p_{*} i^{\flat}$ sends the subcategory $\d{\Sing(Y)}{\coh Y}$
to the subcategory $\d{\Sing(X)}{\coh X}.$ Therefore, we obtain the following corollary.
\begin{corollary}
The functor $\overline{\Phi_{T
*}},$ which realizes an equivalence between the triangulated
categories of singularities of $Y$ and $X,$ gives also a functor
$$
\d{\Sing(Y)}{\coh Y}/\per{\Sing(Y)}{Y}\lto\d{\Sing(X)}{\coh
X}/\per{\Sing(X)}{X},
$$
and this functor is fully faithful.
\end{corollary}
Note that the functor
$
\overline{\Phi_{T
*}}:
\d{\Sing(Y)}{\coh Y}/\per{\Sing(Y)}{Y}\to\d{\Sing(X)}{\coh
X}/\per{\Sing(X)}{X}
$
is not an
equivalence in general.

\section*{Acknowledgments}
I am grateful to  Denis Auroux, Ludmil Katzarkov, Anton Kapustin,
J\'anos Koll\'ar, Amnon Neeman, Tony Pantev,
Leo Alonso Tarr\'io, and Ana Jerem\'ias L\'opez for very useful
discussions.

And, finally, I want to thank my PhD adviser Vasily
Alexeevich Iskovskikh who always supported me  during all my
mathematical life and who passed away on January 4-th, 2009.

\end{document}